\documentclass[12pt]{article}
\usepackage[latin1]{inputenc}
\usepackage{amsthm, amssymb}
\usepackage{amscd}
\usepackage{amsmath}
\usepackage{amsfonts}
\usepackage{mathrsfs}
\usepackage{color}
\newtheorem{theorem}{Theorem}

\newtheorem{lemma}{Lemma}
\newtheorem{corollary}{Corollary}

\newtheorem{remark}{Remark}

\newtheorem{definition}{Definition}

\theoremstyle{definition}
\newtheorem{example}{Example}
\numberwithin{example}{section}
\numberwithin{equation}{section}
\numberwithin{table}{section}
\numberwithin{figure}{section}

\renewcommand{\frac}[2]{\genfrac{}{}{}{}{#1}{#2}}

\renewcommand{\a}{\alpha}

\renewcommand{\L}{\Lambda}

\renewcommand{\hat}{\widehat}

\def\R{\mathbb R}

\def\Pb{\mathbb P}
\def\E{\mathbb E}

\def\LL{{\cal L}}

\def\FF{{\cal F}}

\def\bx{{\boldsymbol x}}

\def\bv{{\boldsymbol v}}

\def\L2{{\mathscr L}^2_\rho(X)}
\def\LL{{\mathscr L}}

\begin{document}
\title{Cucker-Smale Flocking Under Hierarchical Leadership and Random Interactions}
\author{Federico Dalmao\thanks{
Mailing address: Centro de Matem\'atica, Facultad de Ciencias, Igu\'a 4225, 11400, Montevideo, Uruguay.
Telephone number: 598 2 5252522,
E-mail: fdalmao@cmat.edu.uy,
Fax number: 598 2 5220653.} 
\qquad 
Ernesto Mordecki\thanks{
Mailing address: Centro de Matem\'atica, Facultad de Ciencias, Igu\'a 4225, 11400, Montevideo, Uruguay.
Telephone number: 598 2 5252522,
E-mail: mordecki@cmat.edu.uy,
Fax number: 598 2 522 06 53.}
}
\date{}
\maketitle
\begin{abstract}
Consider a flock of birds that fly interacting between them. 
The interactions are modelled through a hierarchical system in which each bird, 
at each time step,
adjusts its own velocity according to his past velocity and a weighted
mean of the relative velocities of its superiors in the hierarchy. 
We consider the additional fact, that each of
the birds can fail to see any of its superiors with certain
probability, that can depend on the distances between them.
For this model with random interactions we prove that the
flocking phenomena, obtained for similar deterministic models,
holds true.
\end{abstract}

\section{Introduction}\label{sec:intro}

Consider a flock of $k$ birds 
with positions and velocities denoted respectively by 
$(\bx_1,\dots,\bx_k):=\bx$ 
and 
$(\bv_1,\dots,\bv_k):=\bv$. 
All individual positions and velocities are 
vectors in $\R^3$, and we refer to $\bx$ and $\bv$ as the position
and velocity of the flock respectively. 

Cucker and Smale, in their seminal paper  \cite{CS06}, propose 
to consider a dynamic for this system in which each bird adjusts 
its velocity according to its own past velocity and a weighted mean of
the relative velocities of the other birds of the flock. More precisely, 
considering a discrete time framework with time step $h>0$,
at time $t+h$ the velocity of
bird $i$ is determined by the velocities of all birds 
at time $t$ by the equation
\begin{equation}\label{eq:model}
   \bv_i(t+h)=\bv_i(t)+h\sum_{j=1}^ka_{ij}(\bv_j(t)-\bv_i(t)),
\end{equation}
for certain weighting coefficients $a_{ij}$, $i,j=1\dots,k$.
The dynamic of the
system is completed by the innovations in the position,
according to the equation
\begin{equation}\label{eq:position}
\bx(t+h)=\bx(t)+h\bv(t).
\end{equation}
To simplify the notation we write $x[t]$ for $x(th)$ and similarly for other functions of $t$.
The interactions that Cucker and Smale \cite{CS06} propose to consider
are of the form
\begin{equation*}
  a_{ij}=a_{ij}[t+1]=\frac{K}{(\sigma^2+\|\bx_i[t]-\bx_j[t]\|^2)^\beta},
\end{equation*}
for positive parameters $K,\sigma,\beta$, where $\|.\|$ is the euclidean norm in $\R^3$. 
We remark that this proposal is related to Vicsek's 
model \cite{VCBC:95}, and refer to \cite{CS06} for a comparison 
between Vicsek and Cucker-Smale models.

A modification of Cucker-Smale model was proposed by Shen \cite{Shen:07},
who considers a non-symmetric structure of interactions, 
where each bird is connected only with its superiors in a prescribed hierarchy. 
More precisely, Shen proposes in \cite{Shen:07} to consider the same coefficients of
the Cucker-Smale model if $j$ is superior to $i$ in the hierarchy and zero otherwise. 
This model was named \emph{Cucker-Smale model under hierarchical leadership}.

It is useful to consider that each bird in the flock constitutes a
vertex of a graph. The interaction between two birds is modelled through
the presence of an edge between the corresponding vertices. 
The hierarchical order imply that the adjacency matrix of the graph can be
taken to be triangular inferior (see Definition \ref{definition:hl} below). Our random selection
of the interactions between birds imply that each entry of the matrix
is a random variable. 

One relevant question for these models is to study the asymptotic
behaviour of the flock, more precisely, to answer under which conditions the 
birds do flock together, 
i.e. under which conditions all individual velocities converge to a 
common velocity and all relative positions of one bird with respect to
another converge to a constant limit. We see that in this situation,
when birds flock, in a coordinate reference
system that travels with the common limit velocity,
each bird tends to a constant position, with asymptotically vanishing velocity.

Our proposal in the present paper, departing from the model proposed by Shen \cite{Shen:07}, 
is to consider that at each time step, the interaction between one bird and
each one of its superior is random, including different possibilities: 
the presence and/or the magnitude of the connectivity can depend on the distance between them.
The precise condition is \eqref{eq:key}
in Theorem \ref{theorem:1}.
See also Examples \ref{example1}-\ref{example2}.
In conclusion, we consider what can be named as the 
\emph{Cucker-Smale model under hierarchical leadership with random interactions}, 
being our purpose to study the flocking phenomena for this model.

Other works including randomness in the flock dynamic are \cite{CM07}, where a random external perturbation is included at each step; and \cite{DM10}, where the connectivity at each step and for each pair of birds has a constant, independent, failure rate (see example \ref{example1}). 

The rest of the paper includes Preliminaries in Section 2, our main result in Section 3, Examples in Section 4,
the Proofs in Section 5, and Conclusions in Section 6.

\section{Preliminaries}

In this section we present our model and some general results for them.

\subsection{Flocks under Hierarchical Leadership (HL-flocks)}

\begin{definition}[HL-flock, \cite{Shen:07}]\label{definition:hl}
We say that a flock of $k$ birds is under 
\emph{hierarchical leadership}, or for short that it is an \emph{HL-flock},
if the birds can be labelled as 
$1,\dots,k$ so that
\begin{enumerate}
\item $a_{ij}\;(=a_{\textrm{bird }i\textrm{ led by bird }j})\neq 0 \textrm{ implies that }j<i$,
\item if we define the leader set of each bird $i$ by $$\LL(i)=\{j:a_{i\,j}\neq 0\},$$
then for any $i>1$, $\LL(i)\neq \emptyset$.
\end{enumerate}
\end{definition}
Note that a flock is an HL-flock if and only if the adjacency matrix $((a_{ij}))$
of the associated graph is \emph{lower triangular} with at least one non
vanishing entry at each row except the first one. 
Besides, bird $1$ is the overall leader of the flock and 
it moves with constant velocity. For $i>1$, 
bird $i$ can look only at birds with a smaller index, and should look at least at one of them.

According to our definition, the velocity of bird $i$ is given by the equation
\begin{equation}\label{eq:hl}
   \bv_i[t+1]=\bv_i[t]+h\sum_{j\in\LL (i)}a_{ij}(\bv_j[t]-\bv_i[t]).
\end{equation}
In order to obtain a weighted mean of the velocities of the
superiors of each bird in the previous formula, 
from now on, and through all the paper, we assume the following two
restrictions on the parameters of the model:
\begin{equation}\label{eq:conditions}
h\leq\frac{1}{k-1}, 
\qquad
0\leq a_{ij}\leq 1. 
\end{equation}
The position of the flock is governed by \eqref{eq:position}.
The model considered in this paper assumes that $\LL(i)$ is fixed for each $i=1.\dots,k$,
and that the coefficients $a_{ij}, j\in\LL(i),i=1,\dots,k$ can be random.

\subsection{Change of coordinates}
In the models considered in \cite{CM07,CS06}, the adjacency matrix $A=((a_{ij}))$
is \emph{symmetric} and
in the proof of flocking a key role is played by the Fiedler number $\phi$,
i.e. the second least eigenvalue of 
the Laplacian of the matrix $A$. 
Under hierarchical leadership, the model we consider,
since the adjacency matrix $A=((a_{ij}))$ 
is lower
triangular, this approach is no longer useful.
We refer to \cite{Shen:07} for details.

To overcome this difficulty we adapt methods from \cite{Shen:07,CD09}.
Observe first that, as $\LL(1)=\emptyset$, it follows that
$\bv_{1}[t]=\bv_{1}[0]$ for all $t$, so, if we have flocking,
all limiting velocities are equal to $\bv_{1}[0]$.  
As a consequence of this fact, we decompose the velocities considering relative
velocities with respect to $\bv_{1}[0]$. 
For each bird $i$ we write its velocity as
\begin{equation*}
\bv_{i}[t]=\left(\bv_{i}[t]-\bv_{1}[0]\right)+\bv_{1}[0],
\end{equation*}
and observe that the relevant dynamic for the properties we are interested in depends only on the differences 
$v_{i}[t]:=\bv_{i}[t]-\bv_{1}[0],\ i=1,\dots,k$. Accordingly, the origin of our coordinate system
will be the first bird, introducing then $x_{i}[t]=\bx_{i}[t]-\bx_{1}[t],\ i=1,\dots,k$.

It is easy to see that the relative velocities $v_{i}$ follow the same dynamic that the absolute velocities $\bv_{i}$. 
In fact, from equation (\ref{eq:hl}) we have
\begin{multline*}
\bv_i[t+1]-\bv_1[0] = \bv_i[t]-\bv_1[0]+h\sum_{j\in\LL(i)}a_{ij}\left(\bv_j[t]-\bv_i[t]\right),\\
=\bv_i[t]-\bv_1[0]+
h\sum_{j\in\LL(i)}a_{ij}\big(\left(\bv_j[t]-\bv_1[0]\right)-\left(\bv_i[t]-\bv_1[0]\right)\big),
\end{multline*}
and similarly for the positions. 
Therefore, we have the following dynamic
\begin{equation}\label{eq:our}
\begin{cases}
  x_i[t+1]=x_i[t]+h v_i[t],\\
  v_i[t+1]=v_i[t]+h\sum_{j\in\LL(i)}a_{ij}\left(v_j[t]-v_i[t]\right).
 \end{cases}
\end{equation}

Despite the fact that $v_1[t]=0$ and $x_1[t]=0$ for all $t$, sometimes 
we include them in the formulas for notational convenience.

\subsection{Basic results}
We now present some basic facts about the general model just described,
in particular, in an HL-flock, 
the norm of the velocities, with an adequate norm, is non increasing
in time. More precisely, for $y=(y_1,\dots,y_k)$ in $\R^{3k}$, 
denote
$$
\left\|y\right\|_{\infty}=\max_{1\leq i\leq k}\left\|y_{i}\right\|.
$$

Denote $x_0=\|x[0]\|_\infty$ and $v_0=\|v[0]\|_\infty$.
\begin{lemma}\label{Dcero}
Consider an HL-flock under the dynamic given in \eqref{eq:hl} 
under the constrains \eqref{eq:conditions}.
Then:
\begin{enumerate}
\item[\rm (a)] 
The absolute velocity of the flock is decreasing in norm, i.e., for all $t\geq
0$,  
we have 
$$
\left\|\bv[t+1]\right\|_{\infty}\leq
\left\|\bv[t]\right\|_{\infty}.
$$
\item[\rm (b)] 
The relative velocity is decreasing in norm,
i.e., for all $t\geq 0$,  we have 
$$
\left\|v[t+1]\right\|_{\infty}\leq
\left\|v[t]\right\|_{\infty}.
$$
\item[\rm (c)] The relative position satisfy
\begin{equation*}
\|x[t]\|_\infty\leq x_0+h v_0 t,
\end{equation*}
\item[\rm (d)] The differences of the relative positions verify
\begin{equation*}
\max_{1\leq i,j\leq k}\|x_i[t]-x_j[t]\|_\infty\leq 2x_0+2hv_0 t.
\end{equation*}
\end{enumerate}
\end{lemma}
The proof of (a) and (b) in this Lemma can be easily adapted from \cite{CD09}. See also Dalmao \cite{FDA}. The proof of (c) and (d) are a consequence of $\|v[t]\|_\infty\leq\|v[0]\|_\infty$ and \eqref{eq:position}.

The norm $\|x\|_\infty$ can be considered as a measure of the dispersion of
the flock and similarly for $\|v\|_\infty$. We then give the following
definition.
\begin{definition} In a dynamical system governed by equations
  \eqref{eq:our}, we have \emph{flocking} when the following two
 conditions are fulfilled:
\begin{itemize}
\item[\rm (a)]  
The relative velocities asymptotically vanishes, i.e.
$$
\|v[t]\|_{\infty}\to 0 \text{ when } t\to\infty.
$$
\item[\rm (b)]  
The relative positions converge, i.e. there exists $\hat{x}\in \R^{3k}$ such that
$$
x[t]\to \hat{x} \text{ when } t\to\infty.
$$
\end{itemize}
\end{definition}
\begin{remark}\label{remark:1} 
It is clear to see that a sufficient condition for flocking is 
\begin{equation}\label{eq:suff}
\sum_{t=0}^{\infty}\|v[t]\|_{\infty}<\infty.
\end{equation}
In case of having a random system, a sufficient condition is 
\begin{equation}\label{eq:suff2}
\sum_{t=0}^{\infty}\E\|v[t]\|_{\infty}<\infty,
\end{equation}
as this condition implies that condition \eqref{eq:suff} holds with probability one.
\end{remark}

\section{Main Result}

A relevant issue in Cucker-Smale proposal is the mathematical formulation of the fact that,
in order to prove flocking, certain control from below on the interactions between birds should be considered. 
With the introduction of randomness in the model arises the question of how to mathematically formulate this fact in an intuitive way being also mathematically tractable. 
Our proposal in \eqref{eq:key} is to express this condition through a conditional expectation on the past of the system, where the vector $\FF_{t}=\big(x[0],v[0],\dots, x[t],v[t]\big)$ contains the information of postions and velocities of the flock up to time $t$.

\begin{theorem}\label{theorem:1}
Consider the HL-flock of Definition \ref{definition:hl} that follows the dynamic \eqref{eq:our},  subject to condition \eqref{eq:conditions} and
\begin{equation}\label{eq:key}
\E(a_{ij}[t+1]\mid\mathcal{F}_t)\geq \frac{p}{\left(1+\|x_i[t]-x_j[t]\|\right)^{\alpha}},
\end{equation}
for $j\in\LL(i),\;i=1,\dots,k$,
and for some $p\in (0,1]$. 
Then, if\\
1) $\alpha<1$, or\\
2) $\alpha=1$ and the coefficients
\begin{equation*}
\gamma_{\ell}:=\sum_{j\in\mathcal{L}({\ell})}\frac{p}{\|v_{\ell}[0]-v_j[0]\|},\quad \ell=2,\dots,k,
\end{equation*}
(where $1/0=\infty$) satisfy the following conditions
\begin{equation}\label{eq:condition}
\gamma_{\ell}>k-\ell+2,\text{ for $\ell=3,\dots,k$,\  and\  }\gamma_2>k-1;
\end{equation}
we have flocking.
\end{theorem}

\begin{corollary}\label{cor:1}
Under the conditions of Theorem \ref{theorem:1}, for $\alpha<1$ and any deterministic sequence $\{\delta_t\}$ such that
\begin{equation}\label{corollary1}
\sum_{t=1}^{\infty} \delta^{-1}_t t^{k-2}
\exp{\left(-\frac{p}{(2v_0)^{\alpha}(1-\alpha)} (ht)^{1-\alpha}\right)}<\infty,
\end{equation}
there exists a random time $\tau_0$ such that
\begin{equation*}
\|v[t]\|_\infty\leq \delta_t,
\end{equation*}
for $t\geq \tau_0$.
\end{corollary}

\begin{corollary}\label{cor:2}
Consider the flock of Theorem \ref{theorem:1}. Then, for $\alpha=1$, under the condition
\begin{equation}\label{corollary2}
v_0<\frac{p}{2(k-1)},
\end{equation}
we have flocking.
\end{corollary}

\section{Examples}
\begin{example}[Random failure of connectivity]\label{example1}
Consider that the presence of connectivity between each pair of birds 
is decided independently of the other birds and of what happened at
previous times, with constant probability for all pair of birds and
times. If the connectivity is present, the strength of this
connectivity, i.e. the coefficient in the weighted mean in \eqref{eq:model}
obeys Cucker-Smale formula. 
The mathematical formulation of this situation
consists in considering
\begin{equation}\label{eq:adj}
  a_{ij}[t+1]=X_{t,ij}\frac{1}{(1+\|\bx_i[t]-\bx_j[t]\|)^\a},
\end{equation}
for some fixed $\alpha> 0$, where the $\{X_{t,ij}\}$ are independent Bernoulli random
variables, with
$$
\Pb(X_{t,ij}=1)=1-\Pb(X_{t,ij}=0)=p,
$$
for some $p\in (0,1)$. It is clear that condition \eqref{eq:key} holds.

This same example is treated by other methods in \cite{DM10}. 
A key role in the approach in \cite{DM10} is played by the law of the iterated logarithm, 
having the advantage of giving a more accurate convergence rate.
\end{example}

\begin{example}[Random interactions]
Assume that $\{X_{t,ij}\}$ in \eqref{eq:adj} are nonnegative, independent random variables that satisfy the condition 
$$\inf_{t,i,j}\E X_{t,ij}\geq p.$$
It is clear that condition \eqref{eq:key} holds.
This is a generalization of example \ref{example1}.
\end{example}

\begin{example}[Connectivity in random environment]\label{example2}
We select the active links for each bird at each
time step considering that the \emph{probability} of connection is a function
of the distance between birds, and the strength of the interaction, in
case of presence, is a constant.
Hence transition probabilities of the system at time $t+1$ 
depend on the vector of all the past values of the positions and velocities $\FF_{t}$. 

More precisely, for an HL-flock we consider coefficients $a_{ij}[t],\,t=0,1,\dots$, with $j\in \LL(i)$, such that they can take the values $0,\,p$ with probabilities
\begin{align*}
\Pb(a_{ij}[t+1]=p\,|\,\FF_t)&=1-\Pb(a_{ij}[t+1]=0|\FF_t)\\
&=\frac{1}{\left(1+\left\|x_{i}[t]-x_{j}[t]\right\|\right)^{\alpha}},
\end{align*}
that verifies condition \eqref{eq:key}.
\end{example}

\section{Proofs}
We begin by an auxiliary result that resumes the contractive nature of the dynamic.
Denote 
$$A_{0}=1+2x_{0},
\qquad 
B_{0}=2hv_{0}, 
\qquad
w_{0\ell}=\min\{\|v_\ell[0]-v_j[0]\|\colon j\in\LL(\ell)\},$$
and assume that $w_{0\ell}>0$ for all $\ell=2,\dots,d$ (see Remark \ref{beta} below).
\begin{lemma}\label{lema:X}
Under the conditions of Theorem \ref{theorem:1}, if $\alpha<1$ we have
\begin{multline*}
\E\left[\prod^{t+1}_{\sigma=\tau+1} \Big(1-h\sum_{j\in\LL(\ell)}a_{\ell j}[\sigma]\Big)\Big|\FF_\tau\right]\\
\leq 
\exp{
\left[-\frac{hp}{(1-\alpha)B_{0}}
\left(\left(A_{0}+B_{0} t\right)^{1-\alpha}-
\left(A_0+B_0 \tau\right)^{1-\alpha}\right)\right]},
\end{multline*}
and in case $\alpha=1$:
\begin{equation*}
\E\left[\prod^{t+1}_{\sigma=\tau+1} \Big(1-h\sum_{j\in\LL(\ell)}a_{\ell j}[\sigma]\Big)\Big|\FF_\tau\right]
\leq 
\left(\frac{A_0+hw_{0\ell} \tau}{A_0+hw_{0\ell}(t+1)}\right)^{\gamma_\ell}.
\end{equation*}
\end{lemma}

\begin{proof} Based on properties of the conditional expectation, we have
\begin{multline}\label{eq:reference} 
\E\left[\prod^{t+1}_{\sigma=\tau+1}\Big(1-h\sum_{j\in\LL(\ell)}a_{\ell j}[\sigma]\Big)\Big|\FF_\tau\right]\\
=
\E\left[\E\Big[\Big(1-h\sum_{j\in\LL(\ell)}a_{\ell j}[t+1]\Big)\prod^t_{\sigma=\tau+1}
\Big(1-h\sum_{j\in\LL(\ell)}a_{\ell j}[\sigma]\Big)\Big|\FF_t\Big]\Big|\FF_\tau\right]\\
=\E\left[\E\Big(1-h\sum_{j\in\LL(\ell)}a_{\ell j}[t+1]\Big|\FF_t\Big)\prod^t_{\sigma=\tau+1}
\Big(1-h\sum_{j\in\LL(\ell)}a_{\ell j}[\sigma]\Big)\Big|\FF_\tau\right].
\end{multline}
Assume that $\alpha<1$. Since $\LL(\ell)\neq\emptyset$ and
$$\E(a_{\ell j}[t+1]|\FF_t)\geq \frac{p}{(1+\|x_j[t]-x_\ell[t]\|)^\alpha}
\geq \frac{p}{(A_{0}+B_{0} t)^\alpha},$$
we have
\begin{equation*}
\E\Big(1-h\sum_{j\in\LL(\ell)}a_{\ell j}[t+1]\Big|\FF_t\Big)
\leq 1-\frac{hp}{(A_0+B_0 t)^\alpha},
\end{equation*}
which is a deterministic bound. Therefore
\begin{multline}\label{espprod}
\E\left[\prod^{t+1}_{\sigma=\tau+1}\Big(1-h\sum_{j\in\LL(\ell)}a_{\ell j}[\sigma]\Big)\Big|\FF_\tau\right]\\
\leq \left(1-\frac{hp}{(A_{0}+B_{0} t)^{\alpha}}\right) 
\E\left[\prod^{t}_{\sigma=\tau+1}\Big(1-h\sum_{j\in\LL(\ell)}a_{\ell j}[\sigma]\Big)\Big|\FF_\tau\right]\\
\leq \prod^{t}_{\sigma=\tau} \left(1-\frac{hp}{(A_{0}+B_{0} \sigma)^{\alpha}}\right)
= \exp\left(\sum^{t}_{\sigma=\tau}\log\left(1-\frac{hp}{(A_{0}+B_{0} \sigma)^\alpha}\right)\right)\\
\leq \exp\left(-hp\sum^{t}_{\sigma=\tau}\frac{1}{(A_{0}+B_{0}\sigma)^\alpha}\right).
\end{multline}
Bounding the sum in the last line of \eqref{espprod} by an integral we obtain the desired bound.

Assume now that $\alpha=1$. Analogously to the case $\alpha<1$ we have
\begin{equation}
\E\left[a_{\ell j}[t+1]\Big|\FF_t\right]\geq \frac{p}{A_{0}+h\|v_{\ell}[0]-v_{j}[0]\|t},
\end{equation}
therefore
\begin{multline}
\E\left[\prod^{t+1}_{\sigma=\tau+1}\Big(1-h\sum_{j\in\LL(\ell)}a_{\ell j}[\sigma]\Big)\Big|\FF_\tau\right]\\
\leq
\exp\left[-h\sum_{j\in\LL(\ell)}\sum^{t}_{\sigma=\tau}
\frac{p}{A_{0}+h\|v_{\ell}[0]-v_{j}[0]\|\sigma}\right]\\
\leq
\exp\left[-p\sum_{j\in\LL(\ell)}
\frac{1}{\|v_{\ell}[0]-v_j[0]\|}
\log\left(\frac{A_{0}+h\|v_{\ell}[0]-v_{j}[0]\|(t+1)}{A_{0}+h\|v_{\ell}[0]-v_{j}[0]\|\tau} \right)\right]\\
\leq \left(\frac{A_{0}+hw_{0\ell}\tau}{A_{0}+hw_{0\ell}(t+1)}\right)^{\gamma_{\ell}}.
\end{multline}
This concludes the proof of the Lemma.
\end{proof}

\begin{proof}[Proof of Theorem \ref{theorem:1}: Case $\alpha<1$ ]
We verify condition (\ref{eq:suff2}) and proceed by induction over $\ell=2,\dots,\,k$, considering the subflocks $\{1,\,2,\dots,\,\ell\}$.

Let us start with $\ell=2$. Recall that the first bird is the overall leader of the flock, 
in particular its (relative) velocity $v_1[t]=0$ for all $t$. 
For the second bird, from Definition \ref{definition:hl} we have $\LL(2)=\{1\}$, then
\begin{align}\label{eq:v2}
v_2[t+1]&=v_2[t]+ha_{21}[t+1](v_1[t]-v_2[t])\notag\\
        &=(1-ha_{2\,1}[t+1])v_2[t]\notag\\
        &= \prod^{t+1}_{\tau=1}(1-ha_{21}[\tau])v_2[0].
\end{align}
By Lemma \ref{lema:X}, we have
\begin{equation*}
\E\|v_2[t+1]\|\leq K\exp\left[
-\frac{hp}{(1-\alpha)B_0}\left(A_0
+B_0 t\right)^{1-\alpha}\right],
\end{equation*}
where $K=\|v_2[0]\|
\exp{\left[-\frac{hp}{(1-\alpha)B_{0}} A_{0}^{1-\alpha}\right]}$, 
this gives condition \eqref{eq:suff2}, concluding the proof for $\ell=2$.

We turn now to the general case. Assume that for $j<\ell$ there exist positive constants $C_j$  such that for all $t$ we have
\begin{equation}\label{ih}
\E\|v_j[t]\|\leq C_jt^{j-2} e^{-\kappa (A_{0}+B_{0}t)^{1-\alpha}},
\end{equation}
where $\kappa=-\frac{hp}{(1-\alpha)B_0}$. 
From \eqref{eq:our}, it follows that
\begin{align}\label{xx}
v_\ell[t+1]&=v_\ell[t]+h\sum_{j\in\LL(\ell)}a_{\ell j}[t+1](v_j[t]-v_\ell[t])\notag\\
&= \left(1-h\sum_{j\in\LL(\ell)}a_{\ell j}[t+1]\right)v_{\ell}[t]+
h\sum_{j\in\LL(\ell)}a_{\ell j}[t+1]v_j[t]\notag\\
&=(1-a[t+1])v_\ell[t]+b[t+1].
\end{align}
where we denote
\begin{equation}\label{ayb}
a[t+1]=\sum_{j\in\LL(\ell)}a_{\ell j}[t+1],\qquad b[t+1]=h\sum_{j\in\LL(\ell)}a_{\ell j}[t+1]v_j[t].
\end{equation}
Iterating \eqref{xx} and taking norms, it follows that
\begin{equation}\label{new}
\|v_\ell[t+1]\|\leq \prod^{t+1}_{\tau=1}(1-ha[\tau])\|v_\ell[0]\|+
\sum^{t+1}_{\tau=1}\|{b}[\tau]\|
\prod^{t+1}_{\sigma=\tau+1}(1-ha[\sigma]),
\end{equation}
where we agree that $\prod^{t+1}_{\sigma=t+2}(1-ha[\sigma])=1$.
By Lemma \ref{lema:X}, there exists a constant $C$ such that
$$\E\prod^{t+1}_{\tau=1}(1-ha[\tau])\leq Ce^{-\kappa (A_{0}+B_{0}t)^{1-\alpha}}.$$
We have
\begin{align*}
\E\|{b}[\tau]\|&\leq h\sum_{j\in\LL(\ell)}\E\left( a_{\ell j}[\tau]\|v_j[\tau-1]\|\right)\\
&\leq h\sum_{j\in\LL(\ell)}\E \|v_j[\tau-1]\|,
\end{align*}
and, by the induction hypothesis \eqref{ih}, for all $t$ we have
\begin{align*}
\sum^{t+1}_{\tau=1} \E\|{b}[\tau]\|
&\leq h\sum_{j\in\LL(\ell)}\sum^{t+1}_{\tau=1}\E \|v_j[\tau-1]\|\\
&\leq h\sum_{j\in\LL(\ell)}\sum^{t+1}_{\tau=1}C_j\tau^{j-2} 
e^{-\kappa(A_{0}+B_{0}(\tau-1))^{1-\alpha}}\\
&\leq C^\prime(t+1)^{\ell-3}\sum^{t+1}_{\tau=1} 
e^{-\kappa(A_{0}+B_{0}(\tau-1))^{1-\alpha}},
\end{align*}
where $C^\prime=\max_{j\in\LL(\ell)}C_j$.
Using this information and Lemma \ref{lema:X}, for the second term on the right hand side of \eqref{new}, we have
\begin{multline*}
\E\left[\sum^{t+1}_{\tau=1}\|{b}[\tau]\|\prod^{t+1}_{\sigma=\tau+1}(1-ha[\sigma])\right]=
\E\left[\sum^{t+1}_{\tau=1}\|{b}[\tau]\|
\E\left[\prod^{t+1}_{\sigma=\tau+1}(1-ha[\sigma])\Big|\FF_\tau\right]\right]\\
\leq C^\prime (t+1)^{\ell-3} \sum^{t+1}_{\tau=1}
e^{-\kappa ((A_{0}+B_{0}t)^{1-\alpha}-(A_{0}+B_{0}\tau)^{1-\alpha})}
e^{-\kappa (A_{0}+B_{0} (\tau-1))^{1-\alpha}}\\
\leq K(t+1)^{\ell-2} e^{-\kappa t^{1-\alpha}},
\end{multline*}
for a convenient constant $K$.
Then, the induction thesis \eqref{ih} for $\ell$ follows with $C_\ell=C+K$. 
This concludes the proof of the Theorem.
\end{proof}

\begin{proof}[Proof of Corollary \ref{cor:1}]
We have
\begin{align*}
\Pb(\|v[t]\|_\infty\geq \delta_t)&\leq \frac{1}{\delta_t}\E\|v[t]\|_\infty
\leq \frac{1}{\delta_t}C_k t^{k-2}e^{-\kappa (A_{0}+B_{0}t^{1-\alpha})},
\end{align*}
and, since the series in \eqref{corollary1} converges, the result follows by the Borel-Cantelli Lemma.
\end{proof}

\begin{proof}[Proof of Theorem \ref{theorem:1}: Case $\alpha=1$] The proof is carried by induction as for the case $\alpha<1$. 
We begin by $\ell=2$.

From \eqref{eq:v2} and Lemma \ref{lema:X} we have
\begin{align*}
\E\|v_2[t+1]\|&\leq
\|v_2[0]\|\E\prod_{\tau=1}^{t+1}\left(1-ha_{21}[\tau]\right)\\
&\leq\|v_2[0]\|\left(\frac{A_{0}}{A_{0}+hw_{0\ell} (t+1)}\right)^{\gamma_2}.
\end{align*}
In view of Remark \ref{remark:1}, 
under the condition $\|v_2[0]\|<p$, the series of the velocities
in \eqref{eq:suff2} is convergent and we have flocking of the first two birds.

Consider now the $\ell$-th bird. For the velocity, we have
\begin{align}
v_{\ell}[t+1]&=v_{\ell}[t]+h\sum_{j\in\mathcal{L}({\ell})}a_{{\ell}j}[t+1]\left(v_j[t]-v_{\ell}[t]\right)\notag\\
        &=\left(1-ha[t+1]\right)v_{\ell}[t]+b[t+1],\label{eq:rew}
\end{align}
where $a[t+1]$ and $b[t+1]$ are defined as in \eqref{ayb}.

With this notation, we have
\begin{equation*}
v_{\ell}[t+1]=\prod_{\tau=1}^{t+1}\left(1-ha[\tau]\right)v_{\ell}[0]
+\sum_{\tau=1}^{t+1}b[\tau]\prod_{\sigma=\tau+1}^{t+1}(1-ha[\sigma]),
\end{equation*}
where we agree that $\prod_{\sigma=t+2}^{t+1}(1-ha[\sigma])=1$.
By Lemma \ref{lema:X} 
we obtain
\begin{align*}
\E\prod_{\sigma=\tau+1}^{t+1}\left(1-ha[\tau]\right)&\leq
\left(\frac{A_{0}+hw_{0\ell}\tau}{A_{0}+hw_{0\ell}(t+1)}\right)^{\gamma_{\ell}}.
\end{align*}
\begin{remark}\label{beta}
The case $w_{0\ell}=0$ is not considered in the previous computations. 
It can be checked with slightly different arguments that in this case the convergence still holds.
\end{remark}

Define recursively $\delta_{\ell}=\delta_{\ell-1}\wedge\gamma_{\ell}-1$.
In what respects $b[t+1]$ in \eqref{eq:rew}, we have
\begin{align*}
\E\|{b}[\tau]\|
&=
h\sum_{j\in\mathcal{L}(\ell)}\E\left(a_{{\ell}j}[\tau]\|v_j[\tau-1]\|\right)
\leq h\sum_{j\in\mathcal{L}(\ell)}\E\|v_j[\tau-1]\|\\
&\leq
h\left(\sum_{j\in\mathcal{L}(\ell)}\|v_j[0]\|\right)
\left(\frac{A_{0}}{A_{0}+hw_{0\ell}\left(\tau-1\right)}\right)^{\delta_{\ell-1}}
\\
&\leq
h\left(\sum_{j\in\mathcal{L}(\ell)}\|v_j[0]\|\right)\left(\frac{A_{0}+hw_{0\ell}}{A_{0}}\right)^{\delta_{\ell-1}}
\left(\frac{A_{0}}{A_{0}+hw_{0\ell}\tau}\right)^{\delta_{\ell-1}}.
\end{align*}
With this information, we conclude that
\begin{multline*}
\E\|v_{\ell}[t+1]\|\leq 
\|v_{\ell}[0]\|
\left(\frac{A_{0}}{A_{0}+hw_{0\ell}(t+1)}\right)^{\gamma_{\ell}}\\
+h\left(\frac{A_{0}+hw_{0\ell}}{A_{0}}\right)^{\delta_{\ell-1}}
\left(\sum_{j\in\mathcal{L}(\ell)}\|v_j[0]\|\right)
(t+1)
\left(\frac{A_{0}}{A_{0}+hw_{0\ell}(t+1)}\right)^{\delta_{\ell-1}\wedge\gamma_{\ell}}\\
\leq k_{\ell}\left(\frac{A_{0}}{A_{0}+hw_{0\ell}(t+1)}\right)^{\delta_{\ell}},
\end{multline*}
where $k_{\ell}$ is a convenient constant that depends on the initial positions and velocities of the flock. If $\delta_{\ell}>1$ we have flocking of the first $\ell$ birds.
It is easy to see that this condition is equivalent to \eqref{eq:condition} with $\ell=k$.
This concludes the proof.
\end{proof}

\subsection{Acknowledgments}
The authors are grateful to Felipe Cucker for helpful discussion and
for providing a preprint of the paper 
``On the critical exponent for flocks under hierarchical leadership'', \cite{CD09}.
\section{Conclusions}
The mathematical modelling of collective behaviour is relevant in many
different scientific frameworks. In our paper we study the existence of
flocking, i.e. the convergence to a common velocity and stable
configuration in the positions of a system of $k$ individuals that
travel, adjusting their own velocities as
a wheighted average of the relative velocities of the other individuals, that, for
simplicity,
we name ``birds''. A very well known model in flocking was proposed by
T.~Vicsek, A.~Czir\'ok, E.~Bet-Jacob, and O.~Shochet in
\cite{VCBC:95},
where each bird is capable of seeing his neighbour birds within a certain fixed
distance to compute the weighted average of velocities. Several
simulation studies were carried out (see for instance \cite{Buhl:06})
but the rigorous proofs of convergence results seems to be
complex, as can be seen in the work by A.~Jadbabaie, J.~Lin, and
A.S. Morse \cite{JLMo:03}. A fundamental advance in the modelling
of the flocking phenomena was the introduction of power decaying (with
the distance)
coefficients for the computation of the weighted average by Cucker and
Smale in \cite{CS06}. This model allows to obtain conditions in terms
of the initial state of the system in order to have flocking. A
modification in the dynamic was proposed by Shen in \cite{Shen:07}
who considers a leadership system in the interactions, in which each
bird has a set of superiors to obey when computing his future velocity
(i.e. a subset of the flock), under the same
Cucker-Smale power decaying interactions. It is clear that the
interactions in Cucker-Smale model are stronger than the ones in the model
proposed in \cite{VCBC:95}. 
Naturally arises then the question, whether it is possible to weaken
the interactions of Cucker-Smale model in order 
to have more realistic modelling, preserving the flocking property.
Our work inscribes in this direction,
incorporating certain probability of failure for each connection at each
time step and random interactions. 
Our main Theorem states that the flocking propoerty 
still holds under random interactions,
provided that the conditional expectation of the interactions remains bounded
from below by the Cucker-Smale interaction coefficient.
It should be noticed that our proposal includes the
possiblity of having at certain times a disconnected flock, but the
probability of this event is handled in such a way that it does not 
prevent the existence (and mathematical treatment) of flocking.

\end{document}